\documentclass{amsart}

\usepackage{amssymb}
 \usepackage{mathtools}
 \usepackage{pifont}
 \usepackage{multicol}
\usepackage{mathrsfs}
\usepackage{multicol}
\usepackage{amsthm}
\usepackage{float}
\usepackage{verbatim} 
\usepackage{tikz}
\usepackage{tikz-cd}
\usetikzlibrary{arrows}
\usetikzlibrary{shapes.arrows}   
\usetikzlibrary{positioning}
\usepackage{mathrsfs}
\usepackage{amssymb}
\usepackage{wasysym}
\usepackage{mathtools}
\usepackage{python}
\usepackage{enumerate}
\usepackage[numbers,comma,sort&compress]{natbib} 
\usepackage[colorlinks=true,citecolor=blue]{hyperref}

\definecolor{1}{rgb}{1,0.2,0.3}
\definecolor{2}{rgb}{0.1,0.3,0.5}
\definecolor{3}{rgb}{1,1,0}
\definecolor{4}{rgb}{255,255,255}
\newcommand{\R}{\mathbb{R}}	   
\newcommand{\N}{\mathbb{N}}
\newcommand{\Z}{\mathbb{Z}}		   
 
\newcommand{\F}{\mathbb{F}}

\newcommand{\one}{\mathbf{1}}
\renewcommand{\Pr}{\mathbb{P}}

\newcommand{\Pro}{\mathbb{P}}
\usepackage{lipsum}

\newtheorem{theorem}{Theorem}[section]
\newtheorem{corollary}[theorem]{Corollary}
\newtheorem{lemma}[theorem]{Lemma}
\newtheorem{conjecture}[theorem]{Conjecture}
\theoremstyle{definition}
\newtheorem{definition}[theorem]{Definition}

\theoremstyle{remark}

\begin{document}

\tikzset
{
  x=1in,
  y=1in,
}

\title{Topology of random $2$-dimensional cubical complexes}

\author{Matthew Kahle}
\address{The Ohio State University}
\email{kahle.70@osu.edu}
\thanks{The first author was supported in part by NSF DMS \#1547357 and CCF \#1740761. He is grateful to the Simons Foundation for a Simons Fellowship, and to the Deutsche Forschungsgemeinschaft (DFG) for a Mercator Fellowship.}

\author{Elliot Paquette}
\address{The Ohio State University}
\email{paquette.30@osu.edu}

\author{\'Erika Rold\'an}
\address{Technische Universit\"{a}t M\"{u}nchen and \`Ecole Polytechnique F\`ed\`erale de Lausanne}
\email{erika.roldan@ma.tum.de}
\thanks{The third author was supported in part by NSF-DMS \#1352386 and NSF-DMS \#1812028. She has also received funding from the European Union’s Horizon 2020 research and innovation program under the Marie Skłodowska-Curie grant agreement No 754462.}



\date{\today}
\maketitle

\begin{abstract}
We study a natural model of random $2$-dimensional cubical complex which is a subcomplex of an $n$-dimensional cube, and where every possible square $2$-face is included independently with probability $p$. Our main result is to exhibit a sharp threshold $p=1/2$ for homology vanishing as $n \to \infty$. This is a $2$-dimensional analogue of the Burtin and Erd\H{o}s--Spencer theorems characterizing the connectivity threshold for random cubical graphs.

Our main result can also be seen as a cubical counterpart to the Linial--Meshulam theorem for random $2$-dimensional simplicial complexes. However, the models exhibit strikingly different behaviors. We show that if $p > 1 - \sqrt{1/2} \approx 0.2929$, then with high probability the fundamental group is a free group with one generator for every maximal $1$-dimensional face. As a corollary, homology vanishing and simple connectivity have the same threshold, even in the strong ``hitting time'' sense. This is in contrast with the simplicial case, where the thresholds are far apart. The proof depends on an iterative algorithm for contracting cycles --- we show that with high probability the algorithm rapidly and dramatically simplifies the fundamental group, converging after only a few steps. 

\end{abstract}

\section{Introduction}
Denote the $n$-dimensional cube by $Q^n = [0,1]^n$, and the set of vertices of the $n$-dimensional cube by $Q^n_0$. In other words, $Q^n_0 = \{0,1\}^n$, which is the set of all $n$-tuples with binary entries. More generally, let $Q^n_k$ denote the $k$-skeleton of $Q^n.$  For example, $Q^n_1$ is the graph with vertex set $Q^n_0$ and an edge (i.e.\ a $1$-dimensional face) between two vertices if and only if they differ by exactly one coordinate.  Define the random $2$-dimensional cubical complex $Q_2(n,p)$ as having $1$-skeleton equal to $Q^n_1$, and where each $2$-dimensional face of $Q^n$ is included with probability $p$, independently. 

The random complex $Q_2(n,p)$ is a cubical analogue of the random simplicial complex $Y_2(n,p)$ introduced by Linial and Meshulam in \cite{LM06}, whose theory is well developed.  The random complex $Y_2(n,p)$ is defined by taking the complete $1$-skeleton of the $n$-dimensional simplex $\Delta^n$ and including each $2$-face independently and with probability $p.$  In this way, $Q_2(n,p)$ is constructed in exactly the same way as $Y_2(n,p),$ \emph{except} that the underlying polytope $\Delta^n$ is replaced by $Q^n.$

The random complex $Q_2(n,p)$ is also the $2$-dimensional analogue of the random cubical graph; see \cite{KSW92} for a 1992 survey. To make the connection precise and review some of the history, let $Q(n,p)$ denote the random subgraph of the $n$-cube defined by including all vertices of $Q^n$, and including each edge in $Q_1^n$ independently with probability $p$.

Burtin \cite{Burtin77}, and later, Erd\H{o}s--Spencer \cite{ES79}, studied the threshold for connectivity in the random cubical graph. 

\begin{theorem}[Burtin \cite{Burtin77} and Erd\H{o}s--Spencer \cite{ES79}]\label{hpr}
For $Q \sim Q(n,p)$ and for any fixed $p \in [0,1]$,
\begin{equation}
    \lim_{n \to \infty}
    \mathbb{P}[Q \text{ is connected}]
    =
    \begin{cases}
        0 &\text{ if } p<1/2,\\[4pt]
        e^{-1} & \text{ if } p=1/2,\\[4pt]
        1 &\text{ if } p>1/2.\\
    \end{cases}
\end{equation}
\end{theorem}
\noindent Let $f_n(p)=\mathbb{P}[Q \sim Q(n,p) \text{ is connected}]$. It is easy to see that if $p < 1/2$, then $\lim_{n \to \infty}f_n(p) =0$. Burtin showed that this is sharp; if $p > 1/2$ then $\lim_{n \to \infty}f_n(p) =1$.

Erd\H{o}s and Spencer refined this argument to show what happens at $p=1/2$ exactly.  They also show that near $p=1/2$, the only connected components of $Q(n,p)$ are isolated vertices and a giant component.  Then it is straightforward to show that the number of isolated vertices has a limiting Poisson distribution with mean $1$, and as a consequence
$$\Pro [{\beta}_0 = k + 1 ] \to e^{-1}/k!$$
for every integer $k \ge 0,$ where ${\beta}_0$ denotes the number of connected components.

This picture strongly parallels what is seen in Erd\H{o}s-R\'enyi graphs $G \sim G(n,p)$. Let $p = (\log n + c)/n$ with $c \in \mathbb{R}$ fixed. Then
\[
   \lim_{n \to \infty}
   \Pr[ G \text{ is connected} ] = e^{-e^{-c}}
\]
(see \cite{ER60} or \cite{Bollobas}).  Letting $c \to \pm \infty$, the probability tends to zero or one. 

The proofs share a strong similarity, in that the method is to enumerate potential cutsets and show that they are rare by making a first moment estimate of the number of cutsets.

It is therefore perhaps reasonable to speculate that the topological phenomenology of the higher dimensional process $Q_2(n,p)$ mirrors that of $Y_2(n,p)$ after appropriately adjusting how $p$ is chosen as a function of $n.$  We shall show, however, that there are major differences between the topology of $Q_2(n,p)$ and $Y_2(n,p).$

Before discussing our results and these differences, we introduce some common terminology. 
Our focus is on typical behavior of random objects for large values of n.
So, we will say that a sequence of statements $\mathcal{P}_n$ holds \emph{with high probability} (abbreviated w.h.p.) if
\[
  \lim_{n \to \infty} \mathbb{P}[\mathcal{P}_n]=1.
\] 

We will make use of the Landau notations $O, o, \omega, \Omega, \Theta$ in the asymptotic sense, so that $f = O(g)$ means $f/g$ is eventually bounded above as $n \to \infty$ and $f = o(g)$ means $f/g \to 0$ as $n \to \infty.$  Also, $f = \omega(g)$ means $g = o(f)$ and $f = \Omega(g)$ means $g = O(f).$ Finally, we will use $f = \Theta(g)$ to mean $f = O(g)$ and $f = \Omega(g)$.  We occasionally display parameters like $O_{a,b,c}(\cdot),$ emphasizing that the implied constants depend on $a,b,c.$

We will also make use of the notion of thresholds. A function $f = f(n)$ is said to be a \emph{threshold} for a property $P$ of a sequence of random objects $G=G_{n,p}$ if $p=\omega(f)$ implies $G \in P$ w.h.p., and $p=o(f)$ implies $G \not\in P$ w.h.p. Such a threshold is only defined up to $n-$independent scalar multiples. If there is a function $g = o(f)$ so that $p \geq f + g$ implies $G \in P$ w.h.p.\ and $p \leq f - g$ implies $G \not\in P$ w.h.p.\ the threshold is \emph{sharp.} If no such $g$ exists, the threshold is coarse.

%
In this paper we study the \emph{fundamental group} $\pi_1(Q)$ for $Q \sim Q_2(n,p).$  The fundamental group can be given a purely combinatorial representation for a space such as $Q_2(n,p),$ which we discuss in Section \ref{sec:pi1}.
Our first result establishes the threshold for $\pi_1(Q)=0$, i.e.\ for $Q\sim Q_2(n,p)$ to be simply connected.  Here, we formulate the theorem for $p$ fixed independently of $n.$  

The $2$-dimensional analogues of isolated vertices are maximal $1$-faces. It is not hard to see that the threshold for maximal $1$-faces is $p=1/2$, and that moreover for any fixed $p \leq 1/2,$ there are maximal $1$-faces w.h.p.  Moreover, for $p=p(n)$ near this threshold they are approximately Poisson distributed. Most of our work is to prove the following.

\begin{theorem}\label{T:vanishing}
If $Q \sim Q_2(n,p)$
and $p>1/2$, then w.h.p.\ $\pi_1 (Q)= 0$.
\end{theorem}

This theorem marks a substantial difference between $Q_2(n,p)$ and $Y_2(n,p).$  The threshold for $\pi_1(Y) = 0$ for $Y \sim Y_2(n,p))$ is roughly $p = n^{-1/2},$ which is proven by Babson, Hoffman, and Kahle \cite{bhk11}.  The threshold is sharpened by Luria and Peled \cite{LuriaPeled18}.  However, for $p=p_n$ below this threshold, also satisfying $p \geq (2+\epsilon) \log n / n$, it is shown in \cite{HKP01} that $\pi_1(Y_2(n,p))$ has Kazhdan's property (T) (see \cite{HKP01} for the discussion therein).  This property precludes the possibility of having any nontrivial free subgroup.

Moreover, the homology vanishing threshold $Q \sim Q_2(n,p)$ \emph{coincides} with the threshold for simple-connectedness.  The first homology group $H_1(Q;\Z)$ is the abelianization $\pi_1(Q)/[ \pi_1(Q),\pi_1(Q)].$  Therefore from Theorem \ref{T:vanishing}, $p = \tfrac 12$ is a sharp threshold for $H_1(Q;\Z) = 0$.

In $Y_2(n,p),$ the homology vanishing threshold is $2\log n / n$ (due to \cite{LM06} over $\F_2$, to \cite{MW09} over general fields, and finally to $\Z$ coefficients by \cite{LP18}; see also \cite{NewmanPaquette01} for the extension to higher dimensions).  Hence in $Y_2(n,p),$ there is a wide range of $p$ \emph{below} the simple-connectivity threshold in which $Y_2(n,p)$ has nontrivial fundamental group but trivial homology.  Again, in $Q_2(n,p)$ the threshold for the vanishing of homology \emph{coincides} with the threshold for the vanishing of $\pi_1.$

For $p > 1-(\tfrac12)^{1/2},$ we are able to completely characterize the fundamental group $\pi_1(Q)$ for $Q \sim Q_2(n,p)):$

\begin{theorem}\label{thm:bubbles}
For $p > 1-(\tfrac12)^{1/2},$ with high probability,
for $Q \sim Q_2(n,p)$
\[
\pi_1(Q) \cong F_N
\]
where $F_N$ denotes a free group on $N$ generators, and $N$ denotes the number of maximal $1$-dimensional faces in $Q$.
\end{theorem}
\noindent An \emph{maximal $1$-face} is a $1$-face which is not contained in any $2$-face.
Hence, just below the $p = \tfrac 12$ threshold, only maximal $1$-faces contribute to the fundamental group.  This is strongly reminiscent of the \emph{homology} of the simplicial complex $Y \sim Y_2(n,p)$ just below its homology vanishing threshold.  From \cite{HKP01}, for $p = t \log n / n$ with $t \in (1, 2),$ the homology group $H_1(Y; \Z)$ is a free abelian group with rank given by the number of isolated $1$-faces of $Y.$  

This also has the following corollary.
\begin{corollary}
If we take for fixed $c \in \R,$ and
\[
p = \frac{1}{2}\biggl( 1 + \frac{\log n + c}{n} \biggr),
\]
then $$\Pr( \beta_1 = k) \to \frac{e^{ck}e^{-e^{-c}}}{k!}.$$
\end{corollary}

Here $\beta_1$ represents the dimension of homology $H_1(Q, \mathbb{R}$. It is also the number of generators of the free group $\pi_1(Q)$. This all follows from the main theorem, together with showing that the number of maximal $1$-faces is Poisson distributed with mean $e^{-c}$.

It is also possible to formulate a \emph{process version} of this statement.  Here one couples all $Q_2(n,p)$ together for all $p \in [0,1]$ in a monotone fashion, so the $2$-faces of $Q_2(n,p_1)$ are a subset of $Q_2(n,p_2)$ whenever $p_1 \leq p_2$.  Let $(Q_p : p \in [0,1])$ have this distribution, and let $N_p$ be the number of maximal $1$-faces in $Q_p.$ Then we can formulate the \emph{stopping times} $T_{sc}$ and $T_{2d}$ as
\[
T_{sc} = \inf\{ p : \pi_1(Q_p) = 0 \}
\quad\text{and}\quad
T_{2d} = \inf\{ p : N_p = 0 \}.\]
We have the following.

\begin{corollary}
$T_{sc} = T_{2d}$ w.h.p.
\end{corollary}

Theorem \ref{thm:bubbles} characterizes the structure of $\pi_1(Q)$ for $Q \sim Q_2(n,p)$ and $$1- \left( \frac{1}{2} \right)^{1/2} < p < \frac{1}{2}.$$ For smaller $p$, we are not able to completely describe the structure of the fundamental group, but we give a partial characterization in terms of its free factorization. 

For a finitely-presented group $H$, we say that $H$ is \emph{indecomposable} if whenever $H$ is written as a free product of two groups $H \cong H_1 * H_2$, either $H_1$ or $H_2$ is a trivial group.  It is well known that a finitely-presented group $G$ be may be written as a free product
$$G=G_1 * G_2 * \dots * G_{\ell},$$
where every $G_i$ is indecomposable. Moreover, this free product factorization is unique, up to isomorphism types of the factors and reordering. 

\begin{definition}\label{def:T_p}
For a cubical sub-complex $T$ of $Q^n,$ let $e(T)$ denote the number of edges in $T.$
Let $\mathscr{T}_p$ be the set of pure 2-dimensional strongly connected cubical complexes $T$ that are subcomplexes of $Q_2^n$ for some $n$ and
so that $1-(1/2)^{1/e(T)} > p.$
\end{definition}

While we do not characterize all the free factors, we are able to characterize some. Essentially, we show that everything that can happen will happen, with high probability.

\begin{theorem}\label{thm:structure}
For $p \in (0,1)$ fixed and $Q \sim Q_2(n,p),$
let the free product factorization of $G = \pi_1(Q)$ be given by
\[
G \cong G_1 * G_2 * \dots * G_\ell.
\]
With high probability, for every $T \in \mathscr{T}_p$ with nontrivial fundamental group, $\pi_1(T)$ appears as a free factor, i.e.\ we have $\pi_1(T) \cong G_j $ for some $1 \le j \le \ell$.
\end{theorem}

Indeed, every finitely presented group $G$ can appear as a free factor (see Section 5 for details).  One can also ask the extremal topological combinatorial question: for given homotopy classes of 2-dimensional complexes $T$, what is the smallest $e(T)$ attainable?  Furthermore, one can ask what is the threshold for specific groups to appear as a free factor in $G.$  In Section 5, we give some examples of specific complexes which we believe to be minimal (we give the torus, the projective plane, and the Klein bottle). 

Using the case of the projective plane, when $0 < p < 1 - (\tfrac 12)^{1/40} \approx 0.017179$ and $Q \sim Q_2(n,p)$, we show that $\pi_1(Q)$ has a $\mathbb{Z}/(2\mathbb{Z})$ free factor with high probability (see Corollary \ref{thm:12}).  This shows in particular that $H_1(Q;\Z)$ has torsion elements for all $p \in (0,p_c),$ where $p_c$ is some critical value in $(0,1).$  The question of whether torsion ever  appears in $H_1(Y; \Z)$ is an open question, although it appears that there is almost always a short burst of enormous torsion \cite{KLNP}. 
\begin{conjecture}
For $p > 1 - (\tfrac 12)^{1/40}$ and $Q \sim Q_2(n,p),$ $\pi_1(Q)$ is torsion free with high probability.
\end{conjecture}
\noindent We further believe it is possible that for $p$ above this threshold, $\pi_1(Q)$ is free with high probability  For the case of $\pi_1(Y)$ with $Y\sim Y_2(n,p)$ the sharp threshold is found by Newman in \cite{N18}, improving on previous work of \cite{CCFK}.

\subsection*{Discussion}

We have not addressed many results about the fundamental group of $\pi_1(Y)$ for $Y\sim Y_2(n,p)$,
which may have interesting analogues in $Q_2(n,p),$ that could further elucidate what appear to be deep differences between simplicial and random cubical complexes.  As the body of literature on  $Y_2(n,p),$ is substantial, we discuss possible directions of interest for questions about $Q_2(n,p).$

Many interesting topological phases of $\pi_1(Q)$ are likely to exist when $p$ tends to $0$ with $n.$  For $Y_2(n,p),$ a particularly rich regime of $p$ is when the mean degree of an edge $np$ tends to a constant. We would expect this regime to be similarly rich for $Q_2(n,p)$ and to name a few transitions that should appear in this regime: the collapsibility threshold \cite{ALLM}, the threshold for a giant shadow \cite{LP16}, and the threshold for $\pi_1(Q)$ to have an irreducible factor in its free product factorization with growing rank.

A natural direction is to consider higher dimensional complexes $Q_d(n,p)$ built in an analogous way to $Y_d(n,p).$   For $Q_d(n,p)$, it may be possible to analyze the higher homotopy groups in a similar fashion to what is done here.

In a different direction, we mention that all of the results we present are about the $n\to \infty$ limit but which also have some content for some large $n$.  These could provide useful results for understanding $2$-dimensional percolation on a sufficiently high dimensional lattice $\Z^n.$  There are some recent related results for such higher dimensional cubical percolation \cite{EL19, HT18, HS16}.

\subsubsection*{Multiparameter generalizations}

In the random cubical graph literature, there is a $2$-parameter model $Q(n; p_0, p_1)$ (see \cite{KSW92} for a survey of some results). First, we take a random induced subgraph of the $n$-cube, where every vertex with probability $p_0$ independently. Then we include each of the remaining edges with probability $p_1$ independently. Bond percolation on the hypercube is the random cubical graph where $p_0=1$ and $p_1$ varies, and site percolation is where $p_1 = 1$ and $p_0$ varies.

It seems natural to form a higher dimensional generalization of this model $Q_2(n; p_0, p_1, p_2)$.
Indeed, Costa and Farber have made a detailed study of the analogous model $Y_2(n ; p_0, p_1, p_2)$ (see \cite{CF17I,CF17II,CF17III}), including many interesting results on the fundamental group.  See also \cite{FN19} wherein new questions about the fundamental group are discussed for this multiparameter model.
Our discussion has been about the special case where $p_0=p_1=1$ and $p_2$ varies.

For, $Q_2(n; p_0, p_1, p_2),$ it is natural to ask if there is a critical surface for homology $H_1$ vanishing in the unit cube.  The case of setting $p_1 = p_2 =1$, and letting $p_0$ vary looks particularly interesting, analogous to the site percolation model. Higher homology is no longer monotone, as in, for example, a random clique complex or Vietoris--Rips complex. Are there separate thresholds for $H_0$ vanishing, $H_1$ appearing, $H_1$ vanishing, and $H_2$ appearing?

\subsection*{Overview and organization}

We begin with Section \ref{sec:prelim} where we define some key notions for working with $Q_2(n,p),$ and we make some elementary estimates about it.
In Section \ref{sec:pi1}, we give a combinatorial definition of $\pi_1.$
In Section \ref{sec:star}, we introduce notation to work with subcomplexes of $Q^n,$ and we introduce the notion the random cubical complex $Q_2(n,p)$ and the random graph $Q(n,p^4).$
In Section \ref{sec:ErdosSpencer}, we summarize some estimates from \cite{ES79} that we need about the random graph $Q(n,p).$
In Section \ref{sec:maximal}, we make estimates for the existence of maximal $1$-faces, and we use this to deduce Theorem \ref{T:vanishing} from Theorem \ref{thm:bubbles}.

In Section \ref{sec:pha}, we introduce an algorithm for identifying contractible $4$-cycles in $\pi_1.$  This algorithm reduces the analysis of $\pi_1$ to determining the topology of small subcomplexes.  In this section, we finish the proof of Theorem \ref{T:vanishing} and give a proof of \ref{thm:bubbles}.  

In Section \ref{sec:generalp}, we show a general structure theorem that describes the
free product factorization of $\pi_1(Q_2(n,p)),$ and we then prove Theorem \ref{thm:structure}. In Section \ref{sec:torsion}, we construct specific complexes which show that certain interesting free factors appear.

\section{Preliminaries}\label{sec:prelim}
\subsection{The edge path group}
\label{sec:pi1}

For subcomplexes $Q$ of $Q^n,$ the fundamental group $\pi_1(Q)$ has a nice combinatorial definition as the \emph{edge path group}, which we now define.

Say that two edges (1-faces) of $Q^n$ are adjacent if they intersect at a vertex.
An \textit{edge-path} in $Q^n$ is defined to be a sequence of edges, for which every consecutive pair are adjacent.

In $Q^n$ any 2-face has a 4-cycle as its boundary.  Conversely, every $4$-cycle in $Q^n$ is the boundary of a $2$-face.  Hence, any two adjacent edges are contained in a unique $2$-face and therefore in a unique $4$-cycle in $Q^n.$

Two edge-paths in $Q$ are said to be \textit{edge-equivalent} if one can be obtained from the other by successively doing one of the following moves: 
\begin{itemize}
  \item[1)] replacing two consecutive adjacent edges by the two opposite edges of the 4-cycle $x$ that contains them, \emph{if} the 2-face that bounds $x$ is in $Q$;
  \item[2)] replacing one edge contained in a 4-cycle $x$ with the other three consecutive edges in $x$, \emph{if} the 2-face that bounds $x$ is in $Q$;
  \item[3)] replacing three consecutive edges in a 4-cycle $x$ with the complementary edge in $x$, \emph{if} the 2-face that bounds $x$ is in $Q$;
  \item[4)] removing an edge that appears twice consecutively or adding an edge that appears twice consecutively.
\end{itemize}

 Define $\overline{0}$ to be the vector with only zero entries in $Q^n_0$. An edge-loop at $\overline{0}$ is an edge-path starting and ending at $\overline{0}$. The random edge-path group $\pi_1(Q)$ is defined as the set of edge-equivalence classes of edge-loops at $\overline{0}$ (with product and inverse defined by concatenation and reversal of edge-loop).

We explore first the extremal cases. If $p=0$ then any $Q \sim Q_2(n,p)$ is equal to $Q^n_1$, that is, with a probability of one the complex $Q$ has not a single 2-face included. Observing that in any graph $G$, the number of independent generators in $\pi_1(G)$ is equal to $E(G)-V(G)+1$, in the case of an element $Q \sim Q_2(n,p)$, we get 
\begin{equation}
    E(C)=2^{n-1}n, \text{ and} \quad V(G)=2^n,
\end{equation}
which implies that the number of independent generators in $\pi_1(Q)$ will be at most $2^{n-1}(n-2)+1$. Thus, when $p=0$ we have that $\pi_1(Q)$ is a free group with $2^{n-1}(n-2)+1$ independent generators, and this is the maximum number of independent generators that the edge-path group of a random 2-cubical complex can attain. This number of independent generators is less than the total number of 4-cycles in $Q^n$ which is $2^{n-3}n(n-1)$. If $p=1$ then any $Q \sim Q_2(n,p)$ will have all the 2-faces included, which implies that $\pi_1(Q)=0$ with a probability of 1. 

\subsection{Star notation and the parallel relation} \label{sec:star}
 In $Q^n$, the four vertices belonging to a 4-cycle have $n-2$ equal entries and two coordinate entries that are not equal in all of them. Denote these non-equal coordinate entries as $i$ and $j$, then we can uniquely represent a 4-cycle using an $n$-tuple with $n-2$ fixed binary values and two $*$. One $*$ will be located on coordinate $i$, and the other will be located on coordinate $j$. As an example, the 4-cycles of $Q^3$ are $\{(0,*,*), (*,1,*), (*,*,1), (*,0,*), (1,*,*), (*,*,0)\},$ with, for instance, $(*,*,0)=\{(1,1,0), (1,0,0), (0,1,0), (0,0,0)\}.$ 
 
For dice, physical realizations of the cube $Q^3$, we have a physical intuitive notion of parallel faces; there are three pairs of parallel faces, and if the die is fair each pair should add up to 7. Using the $*$ notation of 4-cycles introduce above, we extend this notion of parallel faces to $Q^n$. 
 \begin{definition}
 Two 4-cycles in $Q^n$ are \textit{parallel} if they have the two $*$ in the same entries, and if their Hamming distance is 1.
 \end{definition}
 Thus, in $Q^3$ (the cube), there three pairs of parallel 4-cycles: $(0,*,*)$ and $(1,*,*)$, the 4-cycles $(0,*,*)$ and $(1,*,*)$, and the 4-cycles $(0,*,*)$ and $(1,*,*)$. With this parallel notion, we are able to define a binary relation in the set of 4-cycles of a random 2-cubical complex. 
 
 We represent a 3-dimensional cube in $Q^n$ with a vector with $n$ entries, three of which have a fixed $*,$ and the rest of which are binary numbers. If we have two parallel faces $x$ and $y$ that have $*$ in the $i$ and $j$ coordinates and which differ (only) in the binary value of the $k$ coordinate, then the cube that contains them is represented by a vector with the three fixed $*$ in entries $i$, $j$, and $k$, and with the rest of the $n-3$ entries equal to the entries of $x$ (which are also equal to the entries of $y$). 
 
 Observe that the other four 4-cycles of the cube will have all the entries that are not $i$, $j$, or $k$ equal to the entries of $x$, two $*$ in positions either $\{i,k\}$ or $\{j,k\},$ and a binary number in the remaining coordinate. Observe that given two parallel 2-faces in $Q^n$ there is a unique 3-dimensional cube in $Q^n$ that contains them.
 
 \begin{definition}
 Given $Q \sim Q_2(n,p)$, two parallel 4-cycles $x, y \in Q$ are related if the 3-dimensional \textit{cube} that contains them has a 2-face attached to each one of the other four 4-cycles in the cube. We represent this by $x \parallel y$. 
 \end{definition}

 \begin{lemma}\label{L:related}
 If two 4-cycles, $x$ and $y$, are related ($x \parallel y$), then they are edge-equivalent.
 \end{lemma}

\begin{definition}[Graph of parallel related 4-cycles]\label{D:parallel_graph}
Given $Q \sim Q_2(n,p)$, we define its graph of parallel cycles, that we represent by $G[Q]$, as the graph with set of vertices $V^n$ whose elements are all the the 4-cycles in $Q^n$ (there are $2^{n-3}n(n-1)$ 4-cycles), and an edge between two of them if they are related by $\parallel$. 
\end{definition}
Observe that $\parallel$ is reflexive but not transitive. This implies, for instance, that $G[Q^n_2]$ (remember that $Q^n_2$ is the 2-skeleton of $Q^n$) is not the complete graph. We can completely characterize $G[Q^n_2]$.

\begin{lemma}\label{L:partition}
The graph $G[Q^n_2]$ has $\binom{n}{2}$ components, and each one of these components is a $Q^{n-2}_1$ graph.
\end{lemma}

\begin{proof}
Fix an $n>0$. By definition of the relation $\parallel$, a necessary condition for two 4-cycles to be related is to have their two stars in the same position. There are $\binom{n}{2}$ ways of choosing the positions of two stars in a vector of size $n$, which implies that there are at most $\binom{n}{2}$ components in $G[Q^n_2]$. This gives us a partition of the set of vertices that we represent by \[V^n=\bigcup _{i=1}^{\binom{n}{2}} V_i^n.\]

Let $1\leq i\leq {\binom{n}{2}}$, in what follows we prove that the induced subgraph of $V^n_i$ in $G[Q^n_2]$ is isomorphic to $Q^{n-2}_1$. Any element in $V^n_i$ has the two stars in the same position, and the rest of the $n-2$ entries have all possible binary entries. This implies that $| V^n_{i}| = |Q_0^{n-2}|$. Let $\phi: V^n_{i} \to Q_0^{n-2}$ be the natural bijection between these two sets of vertices. It is clear from Definition \ref{L:related} that two 4-cycles  $x$ and $y$ in $V^n_{i}$ are connected in $G[Q^n_2]$ if and only if the Hamming distance of $\phi(x)$ and $\phi(y)$ is 1. This implies that $G[V_i^{n}] \equiv Q^{n-2}_1$.
\end{proof}

For a $Q \sim Q_2(n,p)$ the graph $G[Q]$ is a subgraph of $G[Q^n_2]$. We say that a vertex in $V^n$ is colored if the 4-cycle that it represents has its 2-face present in $Q$, and we say that it is not colored otherwise. 
We use the previous established partition of $V^n$ in $\binom{n}{2}$ sets to denote accordingly the induced subgraphs $G_1, ..., G_{\binom{n}{2}}$ of the graph $G[Q^n_2]$. For a $Q \sim Q_2(n,p)$, this partition defines $\binom{n}{2}$ random subgraphs, that we represent by $G_1[Q], ..., G_{\binom{n}{2}}[Q]$. Then, the edges that are included in $G_i[Q]$, for $1 \leq i\leq {\binom{n}{2}}$ depend on the 2-faces included in $Q$. The next lemma characterizes the probability distribution of each $G_i[Q]$.  

\begin{lemma}\label{L:components}
Let $p^*=p^4$ and $Q \sim Q_2(n,p)$, then, for $1 \leq i\leq {\binom{n}{2}}$, each random graph $G_i[Q]$ is a random graph on $Q^{n-2}_1$, with each edge included independently with probability $p^*$  and with each vertex colored independently with probability p.  
Moreover, the vertex colorings are independent of the edge set.  Using the notation established in the Introduction, the uncolored graph $G_i[Q]$ has the same distribution as $Q(n-2,p^*)$ for all $1 \leq i\leq {\binom{n}{2}}$.
\end{lemma}
\begin{proof}
Let $Q \sim Q_2(n,p)$.  From Lemma \ref{L:partition}, we know that each $G_i[Q]$ is a random graph on $Q^{n-2}_1$. Let $x$ and $y$ be two vertices in $G_i[Q]$ that are connected in $G[Q^n_2]$ and represent this edge by $\overline{xy}$. This implies in particular that $x$ and $y$ are 4-cycles that have, with the previously defined star notation, the $*$ in the same entries and Hamming distance equal to 1. Let $c_{\overline{xy}} \in Q^n$ be the unique 3-dimensional cube that contains $x$ and $y$.  

The probability of $\overline{xy}$ being an edge in $G_i[Q]$ is equal to the probability of the other 4-cycles in $c_{\overline{xy}}$ being covered by 2-faces in $Q$. This event happens with probability $p^*=p^4$ because in $Q_2(n,p)$ each 2-face is added independently with probability $p$. 

Moreover, observe that any of these 4-cycles in $c_{\overline{xy}}$ are not vertices in $G_i[Q]$ because they do not have the two stars in the same location as $x$ (or $y$). This implies the independence between the coloring of the vertices and the inclusion of the edges in $G_i[Q]$. 

Let $\mathcal{C}$ be the set of all 3-dimensional cubes $c_{\overline{xy}}$ with $x$ and $y$ varying among all unordered pairs of vertices in $G_i[Q]$ that are connected in $G[Q^n_2]$. Then, by uniqueness of the cube $c_{\overline{xy}} \in Q^n$, we have that $| \mathcal{C} |$ is equal to the number of edges in $G_i[Q]$. Finally, edges in $G_i[Q]$ are added independently with probability $p^*=p^4$ because each 4-face in a cube in $\mathcal{C}$ only appears in one 3-dimensional cube in $\mathcal{C}$. 
\end{proof}
\subsection{The sizes of the components of $Q(n,p)$} \label{sec:ErdosSpencer}
For a given $p \in (0,1)$, we want to understand the sizes of the components of $Q(n,p)$. We will require an argument from \cite{ES79} that rules out components of small sizes. As we slightly adapt those lemmas, we include proofs below.

Denote by $\mathscr{Q}_s$ the set of all subsets of vertices in $Q^n$ which are connected and have cardinality $s$. Given a subset $S$ of vertices in $Q^n$ that are connected, define
\begin{equation} 
b(S)= \mid\{(u,v)\in Q^n \mid (u,v) \text{ is an edge in } Q^n, \text{ } u \in S, \text{ and } v\notin S  \}\mid.
\end{equation}

Let 
\begin{equation}\label{E:defg}
    g(s)= \sum_{S\in \mathscr{Q}_s}(1-p)^{b(S)}.
\end{equation}
Then from a union bound, $g(s)$ is an upper bound for the probability of the existence of a connected component on $s$ vertices appearing in $Q(n,p).$
\begin{lemma}\label{lem:gs}
\begin{equation}
g(s)\leq 2^n(ns)^s(1-p)^{s(n-\lfloor\log_2(s)\rfloor)}
\end{equation}
\end{lemma}
\begin{proof}
Let $s \geq 1$, then for any $S \in \mathscr{Q}_s$ we have from \cite{hart1976note},
\begin{equation}
b(S) \geq s(n-\lfloor \log2(s)\rfloor ).
\end{equation}
Also, by using that the degree of each vertex of $Q^n$ is at most $n,$
\begin{equation}
    \mid \mathscr{Q}_s \mid 
    \leq 
    2^n (n) (2n) (3n) \cdots ((s-1)n)
    \leq 2^n(ns)^s.
\end{equation}
Hence,
\begin{equation}
    g(s)\leq \sum_{S\in \mathscr{Q}_s}(1-p)^{ s(n-\lfloor \log_2(s)\rfloor )} \leq  2^n(ns)^s (1-p)^{ s(n-\lfloor \log_2(s)\rfloor)}.
\end{equation}
\end{proof}
\begin{lemma}\label{lem:gssum}
For any $p \in (0,1)$, there is a number $T_p \in \mathbb{N}$ and there exists $\delta$, $\epsilon>0$ such that 
\begin{equation}
    \sum_{s} g(s)<2^{-\delta n}
\end{equation}
with the sum overall $s$ so that $T_p\leq s\leq 2^{\epsilon n}$.
\end{lemma}
\begin{proof}
Let $T_p$ be defined by
\[
T_p = \inf_{T \in \N} \{ 2\cdot (1-p)^T \} < 1.
\]

Then for $s \leq 2^{\epsilon n}$ by Lemma \ref{lem:gs},
\[
g(s)
\leq 2^n(ns)^s(1-p)^{s(n-\lfloor\log_2(s)\rfloor)}
\leq 2^n(1-p)^{-2\epsilon n} (n 2^{\epsilon n}(1-p)^{n})^s
\]
for all $n$ sufficiently large.
Then
\[
\begin{aligned}
\sum_{s=T_p}^{\lfloor 2^{\epsilon n}\rfloor} g(s)
&\leq 2^n(1-p)^{-2\epsilon n}
\sum_{s=T_p}^\infty (n 2^{\epsilon n}(1-p)^{n})^s \\
&\leq 2^n(1-p)^{-2\epsilon n} (n 2^{\epsilon n}(1-p)^{n})^{T_p}(1+o(1)),
\end{aligned}
\]
provided $\epsilon$ is chosen so that $2^\epsilon(1-p) <1$ and $n$ is taken large.  By taking $\epsilon$ sufficiently small
\[
\alpha = 2^{1+\epsilon T_p} (1-p)^{T_p-2\epsilon} < 1.
\]
Hence, in terms of $\alpha,$
\[
\sum_{s=T_p}^{\lfloor 2^{\epsilon d}\rfloor} g(s)
\leq \alpha^n n^{T_p} (1+o(1))
\leq 2^{-\delta n}
\]
for some $\delta > 0$ sufficiently small and all $n$ sufficiently large.

\end{proof}

\subsection{The threshold for maximal $1$-faces} \label{sec:maximal}

Any element $Q \sim Q_2(n,p)$ has $2^{n-1}n$ edges, that we represent by \[e_1, e_2, ..., e_{2^{n-1}n},\] with each one of these edges being in $n-1$ different 4-cycles. We represent by $I_{i}$ the indicator function of the event that the edge $e_i$ is maximal, that is, that none of the $(n-1)$ 4-cycles that contain $e_i$ have an attached 2-face. Then, $$\mathbb{E}[I_i]=(1-p)^{n-1}.$$

Let $\mathcal{I}(Q)$ be the random variable that counts the number of maximal edges in $Q$, i.e. $$\mathcal{I}=\sum_{i=1}^{2^{n-1}n}I_i.$$ 
Then 
\begin{equation} \label{E:exp}
\mathbb{E}[\mathcal{I}]=2^{n-1}n(1-p)^{n-1}.
\end{equation}
Observe that if $p=1/2$, then $\mathbb{E}[\mathcal{I}]=n$.

We now prove that Theorem \ref{T:vanishing} follows from Theorem \ref{thm:bubbles}. 
\begin{proof}[Proof of Theorem \ref{T:vanishing}]
We first establish that for $p > \tfrac12,$ $\mathcal{I} = 0$ w.h.p., and for $p \leq \tfrac 12,$ $\mathcal{I} \geq 2$ w.h.p.
For the first claim, the expectation \eqref{E:exp} tends to $0.$  For the second, agin from (\ref{E:exp}), if $(1-p) \geq 1/2$ for a random 2-cubical complex, $Q \sim Q_2(n,p),$
$$ \mathbb{E}[\mathcal{I}]=2^{n-1}d(1-p)^{n-1} \geq n. $$
Thus $\mathbb{E}[\mathcal{I}] \to \infty$ as $n \to \infty$.
Now, we use a second moment argument (see Corollary 4.3.5 of \cite{AS16}) to prove that $\mathbb{P}[\mathcal{I} \geq 2]\to 1$ as $n\to \infty$.

Fix an edge $e_i$. Any other edge $e_j$ such that $I_j$ is not independent from $I_i$, we represent this non-independent relation between edges $e_i$ and $e_j$ by $j \sim i$, will be an edge of one and only one of the $(n-1)$ 4-cycles that contain $e_i$. There are $3(n-1)$ such edges and $\mathbb{P}[I_j \mid I_i] = (1-p)^{n-2}$.  If we define $\Delta^*_i=\sum_{j\sim i} \mathbb{P}[I_j \mid I_i]$, then
$$\Delta^*_i=\sum_{j \sim i} \mathbb{P}[I_j\mid I_i]=3(n-1)(1-p)^{n-2}.$$
Thus $\Delta^*_i=o(\mathbb{E}[I_n])$ which implies that $\mathbb{P}[\mathcal{I}>n/2] \to 1$ as $n \to \infty.$ 

Hence from Theorem \ref{thm:bubbles}, we have that for $p > \tfrac 12,$ $\pi_1(Q_2(n,p)=0$ with high probability and for any $p<\tfrac 12,$ $\pi_1(Q_2(n,p)) = G * \Z * \Z$ for some group $G$ with high probability.  The event that $Q_2(n,p)$ has such a free factorization is  a \emph{decreasing} event, in that for any complex $Q$ that satisfies $\pi_1(Q) = G * \Z * \Z$ for some group $G,$ removing any $2$-face (i.e.\ removing relations from $\pi_1(Q)$) yields a complex $Q'$ so that $\pi_1(Q') = G' * \Z * \Z$ for some other group $G'.$  It follows that for any $p \leq \tfrac 12,$  
\[
\Pr[\exists~G: \pi_1(Q_2(n,p)) = G * \Z * \Z]
\geq
\Pr[\exists~G: \pi_1(Q_2(n,\tfrac12)) = G * \Z * \Z]
\to 1,
\]
as $n \to \infty,$ which completes the proof.
\end{proof}





\section{Parallel homotopy algorithm}
\label{sec:pha}
In this section, we introduce a simple iterative algorithm for finding contractible 4-cycles.  For $Q_2(n,p)$ with $p > 0,$ this algorithm rapidly and dramatically simplifies the fundamental group to its nontrivial parts.

We begin by introducing the algorithm.  We have defined $V^n$ as the set of all 4-cycles in $Q^n_1$. For any subset $V \subset V^n$, we define the \emph{graph of parallel related 4-cycles} denoted by $G(V)$ in a similar fashion to \ref{D:parallel_graph}: the vertex set of $G(V)$ is given by the $V^n$ and two $4$-cycles $x$ and $y$ are connected if they have stars in the same positions, are contained in a $3$-cube $c,$ and all other $4$-cycles in $c$ are in $V$.

Given a $Q \sim Q_2(n,p)$, we denote by $V^n_1$ the subset of $V^n$ which contains the boundaries of 2-faces in $Q$. We then iteratively run the following procedure, with $t \in \mathbb{N}.$

\noindent \textbf{Stage t:} Build the graph of parallel related 4-cycles $G(V^n_t)$. Define the set of 4-cycles $V^n_{t+1}$ as the set of 4-cycles that are connected in $G(V^n_t)$ to a 4-cycle that is in $V^n_t$.

\noindent The algorithm stops at the first $t$ for which  $V^n_{t+1}=V^n_t$.

As an aside, we observe that half of Theorem \ref{T:vanishing} follows from the following result:
\begin{theorem} \label{thm:stage2}
For $p>1/2$ 
\begin{equation}
    \lim_{n \to \infty}\mathbb{P}[V^n_3=V^n]=1.
\end{equation}
\end{theorem}

\subsection{Stage 1: explosive growth}

For any set of $4$-cycles $V \subset V^n$, say that a set of vertices $S$ in $G(V)$ is a \emph{quasicomponent} if $S$ is connected in $G(V^n)$ and $S$ is disconnected from its complement in $G(V).$
\begin{theorem} \label{thm:stage1}
Let $Q \sim Q_2(n,p)$, and let $A_s$ be the event that a quasicomponent of size $s$ in $G(V^n_1)$ exists. Then for any $p \in (0,1),$ there is an integer $T_p$ and $\epsilon,\delta >0$ so that for all $n$ sufficiently large
\begin{equation}\label{eq:noquasi}
\mathbb{P}\bigl[ \cup_{s=T_p}^{2^{\epsilon n}}A_s \bigr]<2^{-\delta n}.
\end{equation}
Also, the probability that there exists a component of $G(V_1^n)$ bigger than $T_p$ with no vertex in $V^n_1$ tends to zero with $n$. 
\end{theorem}

\begin{proof}
The first part of the statement, inequality \eqref{eq:noquasi}, follows by a union bound and Lemma \ref{lem:gssum} by observing that $\mathbb{P}[A_s]\leq g(s)$ --See equation \eqref{E:defg}.

By virtue of \eqref{eq:noquasi}, it remains to show that there are no components of $G(V_1^n)$ bigger than $2^{\epsilon n}$ that do not intersect $V_1^n.$ Let $W$ be the event that there exists a component in $G(V_1^n)$ of size bigger or equal than $2^{\epsilon n}$ that does not intersect $V_1^n$. We show in what follows that $\mathbb{P}[W]\to 0$ as $n \to \infty$. First, we observe that $G(V_1^n)=G[Q]$ and that the vertices in $V_1^n$ are precisely the colored vertices in $G[Q],$ which by Lemma \ref{L:components} are colored independently with probability equal to $p$. For $1 \leq i \leq {\binom{n}{2}}$, define $W_i$ as the event that there exists in $G_i[Q]$ a component of size bigger or equal to $2^{\epsilon n}$ that has all its vertices uncolored. Thus, by Lemma \ref{L:components},
\begin{equation} \label{E:W}
W = \bigcup_{i=1}^{\binom{n}{2}} W_i.
\end{equation}
 Let $1 \leq i \leq {\binom{n}{2}}.$ Conditioned on knowing $G_i[Q]$, in particular on knowing that there are exactly $l$ components with uncolored vertices and with sizes $s_1, s_2, ..., s_l,$ bigger than $2^{\epsilon n}$ in $G_i[Q]$  we get 
 \begin{equation}\label{E:components}
 \mathbb{P}[W_i \mid G_i] \leq \sum_{k=1}^{l}(1-p)^{s_j}.
 \end{equation}
Observing that $s_1 + s_2 + \cdot \cdot \cdot s_l \leq 2^{n-2}$, it has to be the case that $l\leq 2^{n-2}$, and because $(1-p)< 1$ we have that $(1-p)^{s_{k}} \leq (1-p)^{2^{\epsilon n}}$ for all $1\leq k \leq l$. Thus, from equation \eqref{E:components} we get    
\begin{equation}
    \mathbb{P}[W_i \mid G_i] \leq \sum_{k=1}^{l}(1-p)^{2^{\epsilon n}} \leq 2^{n-2}(1-p)^{2^{\epsilon n}}. 
\end{equation}
This implies that $ \mathbb{E}[\mathbb{P}[W_i \mid G_i]] \leq 2^{n-2}(1-p)^{2^{\epsilon n}}$, and thus 
\begin{equation} \label{E:probabilitiw_1}
\mathbb{P}[W_i]\leq 2^{n-2}(1-p)^{2^{\epsilon n}}
\end{equation}
for all $1 \leq i \leq {\binom{n}{2}}$.  
Finally, by a union bound argument on \eqref{E:W} and inequality \eqref{E:probabilitiw_1} we have that 
\begin{equation}
\mathbb{P}[W] \leq {\binom{n}{2}} 2^{n-2}(1-p)^{2^{\epsilon n}}, 
\end{equation}
with
\begin{equation}
    \lim_{n \to \infty} {\binom{n}{2}} 2^{n-2}(1-p)^{2^{\epsilon n}}= 0.
\end{equation}
\end{proof}

\subsection{Stage 2: Only local defects remain}
Let $\mathcal{F}$ be the event that there is no quasicomponent of $G(V_1^n)$ bigger than $T_p$ that is disjoint from $V_1^n.$  This event was shown to hold w.h.p.\ by Theorem \ref{thm:stage1}.  
\begin{lemma}\label{lem:F}
In the event $\mathcal{F}$, any $4$-cycle $v$ with at least $T_p$ neighbors in $G(V_2^n)$ is in $V_3^n.$  Likewise, any $4$-cycle $v$ with at least $T_p$ neighbors in $G(V_1^n)$ is in $V_2^n.$
\end{lemma}
\begin{proof}
Suppose $\mathcal{F}$ holds, and let $v$ be any $4$-cycle.
Suppose that $v$ has at least $T_p$ neighbors in $G(V_1^n).$  
Then the connected component of $v$ in $G(V_1^n)$ has at least $T_p$ neighbors, and therefore this connected component intersects $V_1^n.$
It follows by the definition of $V_2^n$ that $v \in V_2^n.$

Suppose now that $v$ has at least $T_p$ neighbors in $G(V_2^n).$ 
We may suppose that $v$ is not in a component of $G(V_1^n)$ that intersects $V_1^n,$ for if it were, then $v \in V_2^n$ and we are done.
If none of these neighbors are in $V_2^n,$ then each is in a component of $G(V_1^n)$ disjoint from $V_1^n.$  Hence, the union of these components and the component of $G(V_1^n)$ containing $v$ is a quasicomponent of $G(V_1^n)$ that is disjoint from $V_1^n.$  Moreover, it is a quasicomponent which is larger than $T_p,$ which is disjoint from $V_1^n.$  This does not exist in $\mathcal{F},$ and therefore $v$ has a neighbor in $V_2^n.$  Hence $v \in V_3^n.$
\end{proof}

We will show that as a consequence of Lemma \ref{lem:F}, in Stage 2, all those 4-cycles whose every constituent edge has a high enough degree will be collapsed.  For any $p,$ define 
\begin{equation} \label{eq:lightedge}
M_p
=\inf_{M > 0}
\Pr ( \text{Binomial}( \lfloor M/4 \rfloor, p^3) < T_p) < (\tfrac{1}{2})^{1/4}.
\end{equation}
For any $1$-face $f$ in $Q^n,$ define $\deg(f)$ as the number of $2$-faces in $Q$ containing $f.$
Call a 1-face of $Q \sim Q_2(n,p)$ \emph{light} if its degree is less than or equal to $M_p.$ Otherwise, call it \emph{heavy}.  We show that 4-cycles made from heavy edges are contracted in the second stage of the algorithm:
\begin{lemma}\label{lem:stage2}
For any $p \in (0,1),$ with probability tending to $1$ as $n\to \infty,$
every $4$-cycle whose every 1-face is heavy is contained in $V_3^n.$
\end{lemma}

We will introduce some additional notation for working with faces of $Q$.  For two disjoint sets $U, W \subset [n],$ let $(U^*, W^1)$ denote the $|U|$--dimensional face of $Q$ with $*$s in the positions given by $U,$ and $1$s exactly in the positions given by $W.$

Using symmetry it will be enough to analyze the $4$-cycle $(\{1,2\}^*, \emptyset^1).$  With the $M_p$ from \eqref{eq:lightedge}, define $\mathcal{E}$ as the event that all the 1-faces in the $4$-cycle $(\{1,2\}^*, \emptyset^1)$ are heavy, i.e.
\[
\begin{aligned}
\mathcal{E} = 
\bigl\{
&\deg( (\{1\}^*, \emptyset^1)) > M_p,
\deg( (\{1\}^*, \{2\}^1)) > M_p, \\
&\deg( (\{2\}^*, \emptyset^1)) > M_p,
\deg( (\{2\}^*, \{1\}^1)) > M_p
\bigr\}.
\end{aligned}
\]

To prove Lemma \ref{lem:stage2}, it suffices to show that
\begin{lemma}
For any $p \in (0,1),$ there is an $\epsilon > 0$ so that
\[
\mathbb{P}[\mathcal{E} \cap \mathcal{F} \cap \{ \text{the degree of $( \{1,2\}^*,\emptyset^1)$ in $G(V_2^n)$ is less than $T_p$}\}] \leq n^{O(1)}2^{-(1+\epsilon)n}.
\]
\end{lemma}

\begin{proof} 
The possible neighbors of $( \{1,2\}^*,\emptyset^1)$ in $G(V^n)$ all have the form $( \{1,2\}^*,\{j\}^1)$ for some $3 \leq j \leq n.$
To have an edge between these $4$-cycles in $G(V^n_2),$ we must have that
\[
\begin{aligned}
&( \{1,j\}^*, \emptyset^1) \in V^n_2, \quad  
( \{1,j\}^*, \{2\}^1) \in V^n_2, \\
&( \{2,j\}^*, \emptyset^1) \in V^n_2, \quad 
( \{2,j\}^*, \{1\}^1) \in V^n_2.
\end{aligned}
\]
On the event $\mathcal{F},$ we must only lower bound the degree of these $4$-cycles in $G(V^n_1)$ to ensure they are in $V^n_2.$  Hence, define 
\begin{equation} \label{eq:Y}
\begin{aligned}
    &Y_{1j} = \one\{ \deg(  (\{1,j\}^*, \emptyset^1) ) \geq T_p \}, \quad 
    Y_{2j} = \one\{\deg(  (\{1,j\}^*, \{2\}^1) ) \geq T_p \}, \\
    &Y_{3j} = \one\{\deg(  (\{2,j\}^*, \emptyset^1) ) \geq T_p\}, \quad 
    Y_{4j} = \one\{\deg(  (\{2,j\}^*, \{1\}^1) ) \geq T_p\}.
\end{aligned}
\end{equation}
The $\deg$ above refers to the degree of the $4$-cycle in $G(V^n_1).$
We would like to show there are at least $T_p$ choices $j$ for which all $Y_{\ell j}$ for $\ell \in \{1,2,3,4\}$ are $1.$

 On the event $\mathcal{E},$ there are 4 disjoint sets $R_\ell \subset \{3,4, \dots, d\}$ for $\ell \in \{1,2,3,4\}$ of size $\lfloor M_p/4 \rfloor$ so that
 \[
 \begin{aligned}
  &( \{1,k\}^*, \emptyset^1) \in V^n_1, \text{ for } j \in R_1,
  \quad
  ( \{1,k\}^*, \{2\}^1) \in V^n_1, \text{ for } j \in R_2, \\
  &( \{2,k\}^*, \emptyset^1) \in V^n_1, \text{ for } j \in R_3,
  \quad 
  ( \{2,k\}^*, \{1\}^1) \in V^n_1, \text{ for } j \in R_4.
 \end{aligned}
 \]
 
 Observe that the possible neighbors of $(\{1,j\}^*, \emptyset^1)$, for $j \in \{3,4, \dots, n\}$ are given by $(\{1,j\}^*, \{k\}^1)$ for $k \not\in \{1,j\}.$  For simplicity, we will also discard the case $k=2.$ To have this edge in $G(V^n_1),$ we would need that 
 \[
\begin{aligned}
&( \{1,k\}^*, \emptyset^1) \in V^n_1, \quad  
( \{1,k\}^*, \{j\}^1) \in V^n_1, \\
&( \{j,k\}^*, \emptyset^1) \in V^n_1, \quad 
( \{j,k\}^*, \{1\}^1) \in V^n_1.
\end{aligned}
\]
In particular, for $k \in R_1,$ the first of these requirements is guaranteed.  Hence we can define
\[
Z_{1jk} = \one\{
( \{1,k\}^*, \{j\}^1) \in V^n_1,
( \{j,k\}^*, \emptyset^1) \in V^n_1,
( \{j,k\}^*, \{1\}^1) \in V^n_1
\},
\]
and define 
\[
Z_{1j} = \sum_{k \in R_1} Z_{1jk}.
\]
Then $Z_{1j}$ is a lower bound for $\deg(  (\{1,j\}^*, \emptyset^1) ),$ and so if $Z_{1j}$ is at least $T_p,$ then $Y_{1j}=1.$

We do a similar construction for $\ell \in \{2,3,4\},$ making appropriate modifications.  We list these for clarity below:
\[
\begin{aligned}
&Z_{2jk} = \one\{
( \{1,k\}^*, \{2,j\}^1) \in V^n_1,
( \{j,k\}^*, \{2\}^1) \in V^n_1,
( \{j,k\}^*, \{1,2\}^1) \in V^n_1
\}, \\
&Z_{3jk} = \one\{
( \{2,k\}^*, \{j\}^1) \in V^n_1,
( \{j,k\}^*, \emptyset^1) \in V^n_1,
( \{j,k\}^*, \{1\}^1) \in V^n_1
\}, \\
&Z_{4jk} = \one\{
( \{2,k\}^*, \{1,j\}^1) \in V^n_1,
( \{j,k\}^*, \{1\}^1) \in V^n_1,
( \{j,k\}^*, \{1,2\}^1) \in V^n_1
\}. \\
\end{aligned}
\]
In terms of these, we set $Z_{\ell j} = \sum_{k \in R_\ell } Z_{\ell jk}.$

Let $J = \{3,4, \dots, d\} \setminus (\cup_1^4 R_\ell).$  Then the family
\[
\{ Z_{\ell jk} : \ell \in \{1,2,3,4\}, j \in J, k \in R_\ell \}
\]
are independent random variables.
Moreover for any $\ell \in \{1,2,3,4\}$ and $j \in J,$ from \eqref{eq:lightedge},
\[
\Pr(Z_{\ell j} < T_p)
<
\Pr ( \text{Binomial}( \lfloor M_p/4 \rfloor, p^3) < T_p) 
\leq
(\tfrac{1}{2})^{1/4}.
\]
It follows that with
\[
Z = \sum_{j \in J} \prod_{\ell = 1}^4 \one\{ Z_{\ell j} \geq T_p\},
\]
and with $q = \Pr(Z_{\ell j} \geq T_p)^4 > \tfrac12,$
\[
\Pr(Z < T_p) \leq 
\Pr( \text{Binomial}( n - 3 - M_p, q) < T_p)
=n^{O(1)}(1-q)^{n},
\]
which completes the proof.
\end{proof}

\subsection{Stage 3: The final squeeze}
\label{sec:stage3}
In this section, we draw conclusions on what remains non-contracted in the complex in the third stage.

\subsubsection{The simply connected regime, $p > \tfrac12$}
We begin by showing that for $p > 1/2,$ there are simply no light $1$-faces.  Hence in fact for $p > \tfrac12,$ $V_3^n=V_n$ with high probability (proving Theorem \ref{thm:stage2}).
\begin{lemma} \label{lem:edgedegree}
For any $p > 1/2,$ there is an $\epsilon > 0$ so that with probability tending to $1$ with $n,$ for every $1$-face $f$ of $Q \sim Q_2(n,p),$ $\deg(f)>M_p.$
\end{lemma}
\begin{proof}
The degree of a $1$-face is distributed as $\text{Binomial}(n-2,p)$.  For $p > \tfrac12,$ the probability this is less than any fixed constant $M$ is $n^{O(1)}(1-p)^{n}$.  Hence by a union bound, the lemma follows. 
\end{proof}

\subsubsection{Completely shielded $1$-faces}
Call a 1-face $f \in Q \sim Q_2(n,p)$ \emph{completely shielded} if every 3-face $c \in Q^n$ that contains $f$ only contains heavy $1$-faces of $Q$, besides possibly $f$.  Completely shielded 1-faces modify the fundamental group of $Q$ in a simple way, contributing a free factor of $\Z$ if $f$ is maximal.

To see this we begin with the following definition:
\begin{definition}\label{def:bubble}
Let $f$ be any $1$-face of $Q^n.$
Define the \emph{$n$-bubble around $f$} to be the subcubical complex of $Q^n$ given by 
the union of the complete $1$-skeletons of all $3$-faces containing $f,$ and every $2$-face on this skeleton which does \emph{not} contain $f$.
\end{definition}
\noindent A $n$-bubble has fundamental group $\Z.$
\begin{lemma}\label{lem:dbubble}
For any $n \geq 3,$ and any $n$-bubble $X$ around $f,$
\[
\pi_1(X) \cong \Z.
\]
Furthermore, the complex $X \setminus \{f\}$ and the complex $X \cup \{e\},$ where $e$ is any $2$-face containing $f$, are simply connected.
\end{lemma}

\begin{proof}
Without loss of generality, suppose that $f$ is the face $(\{1\}^*,\emptyset^1)$.  The $3$-faces containing $f$ all have the form $(\{1,i,j\}^*, \emptyset^1),$ and so the $1$-skeleton of $X$ is 
\[
  \{ (\{ i \}^*, A^1)
  : 
  A \subset \{1,2,\ldots,n\},
  i \not\in A, |A \cup \{i\}| \leq 3
  \}.
\]
We claim that all the $4$-cycles containing $f$ are homotopic.  As all other $4$-cycles are contractible from the definition of $X$, the statements in the lemma follow.

The $4$-cycles that contain $f$ are boundaries of the $2$-faces of $Q^n$ of the form
\[
\{ (\{ 1,i \}^*, \emptyset^1)
  : 
  2 \leq i \leq n
  \}.
\]
For any $2 \leq i < j \leq n$, the $3$-face $c=(\{1,i,j\}^*, \emptyset)$  intersected with $X$ contains $4$ $2$-faces.  Moreover, the $2$-faces $(\{1,i\}^*, \emptyset)$ and $(\{1,j\}^*, \emptyset)$ are adjacent in this cube.  Hence, these cycles can be deformed through $c$ to one another.  As this held for any such $i$ and $j,$ the proof follows.
\end{proof}

\begin{lemma}\label{lem:completelyshielded}
For any $p \in (0,1),$
let $\hat Q$ be the cubical complex that results from deleting from $Q \sim Q_2(n,p)$ every completely shielded $1$-face $f$ and any $2$-face of $Q$ containing $f$.  Then with high probability,
\[
\pi_1(Q) \cong \pi_1(\hat Q) 
*
\underbrace{( \Z * \Z * \cdots * \Z )}_{N}
\]
where $N$ denotes the number of completely shielded $1$-faces in $Q$ that are isolated.
\end{lemma}
\begin{proof}
From Lemma \ref{lem:stage2}, all $4$-cycles whose every 1-face is heavy are contractible.  In particular, we do not modify the fundamental group of $Q$ if we include all those $2$-faces into $Q$ whose boundary is in $V_3^n$.  Let $\hat Q$ be this cube complex.  

We now remove completely shielded $1$-faces from $\hat Q$ one at a time, tracking the changes to the fundamental group.  We will show what happens after removing the first.  It will be clear that by using induction, a similar analysis would give the claim in the lemma.

Let $f$ be a completely shielded $1$-face of $\hat Q.$
Let $Q_1$ be the complex that results after removing $f$ from $\hat Q$ and any $2$-face containing $f$.  Let $Q_2$ be the union of all the complete $2$-skeletons of all $3$-faces that contain $f$.  Then $Q_2$ contains a $n$-bubble, and it is exactly a $n$-bubble if $f$ is isolated.

As $Q_2 \cup Q_1 = \hat Q$ and $Q_1 \cap Q_2$ is open and path connected (c.f.\ Lemma \ref{lem:dbubble}, as this complex is a $n$-bubble with its central $1$-face deleted).  Moreover, every $4$-cycle in $Q_1 \cap Q_2$ is contractible, and so $\pi_1(Q_1 \cap Q_2)$ is trivial.  From the Siefert-van Kampen theorem, we therefore have that
\[
\pi_1(\hat Q) 
\cong \pi_1(Q_1) * \pi_1(Q_2).
\]
If $f$ is maximal then from Lemma \ref{lem:dbubble}, the fundamental group $\pi_1(Q_2)$ is isomorphic $\Z.$
\end{proof}

\subsubsection{The velvety bubble phase}
For $p > 1-(\tfrac12)^{1/2} \approx 0.292893,$
we further show that the fundamental group completely reduces to its maximal 1-faces.  In this phase, while light 1-faces may exist in $Q_2(n,p)$ (for $p \leq \tfrac12$), they are well separated.

\begin{lemma}\label{lem:nocloselightweights}
For $p > 1-(\tfrac12)^{1/2},$ with high probability, there are no 3-faces $c \in Q^n$ that contain more than one light 1-face of $Q \sim Q_2(n,p)$.
\end{lemma}
\begin{proof}
For a fixed $c$ and a fixed choice of two 1-faces $f_1,f_2,$ for the degrees of $f_1$ and $f_2$ are both light with probability at most $n^{O(1)}(1-p)^{2n-1}$.  Hence for any $p$ as in the statement of the lemma, there is an $\epsilon > 0$ so that the probability this occurs is $2^{-(1+\epsilon)n + O(\log n)}$. As there are $2^n n^{O(1)}$ many ways to pick a $3$-face with two designated edges, the lemma follows from a first moment estimate. 
\end{proof}

We now give the proof of Theorem \ref{thm:bubbles}, which we recall for convenience.
\begin{theorem}
For $p > 1-(\tfrac12)^{1/2},$ with high probability,
for $Q \sim Q_2(n,p)$
\[
\pi_1(Q) \cong
\underbrace{( \Z * \Z * \cdots * \Z )}_{N},
\]
where $N$ denotes the number of maximal $1$-faces in $Q.$
\end{theorem}
\begin{proof}

From Lemma \ref{lem:stage2}, with high probability every $4$-cycle containing only heavy $1$-faces is in $V_3^n$.   From Lemma \ref{lem:nocloselightweights}, with high probability no 3-faces $c \in Q^n$ contain more than one light face.  Hence taking $\tilde Q$ as $Q$ together with all $2$-faces bounded by some element of $V_3^n$ (so that $\pi_1(\tilde Q) = \pi_1(Q)$) every light 1-face $f$ of $\tilde Q$ is completely shielded in $\tilde Q.$  Moreover, every $4$-cycle of $\tilde Q$ either intersects a light $1$-face, or it is the boundary of a 2-face.  Hence in the notation of Lemma \ref{lem:completelyshielded}, $\pi_1(\hat Q)=0$.  It follows that from Lemma \ref{lem:completelyshielded} is a free group on $N'$ generators, with $N'$ the number of completely shielded maximal $1$-faces.  As every light $1$-face is completely shielded w.h.p, it follows that $N'=N$ with high probability.
\end{proof}
 
\section{Structure theorem for general $p$}
\label{sec:generalp}

In this section, we prove Theorem \ref{thm:structure}.
For convenience, we recall some definitions from the introduction. Recall Definition \ref{def:T_p}:
\begin{definition}
For a cubical sub-complex $T$ of any cube $Q^n,$ let $e(T)$ denote the number of edges in $T.$
Let $\mathscr{T}_p$ be the set of pure 2-dimensional strongly connected cubical complexes $T$ that are subcomplexes of $Q_2^n$ for some $n$ and
so that $(1-(\frac12)^{1/e(T)}) < p.$
\end{definition}
\noindent We will prove Theorem \ref{thm:structure}, which we recall below:
\begin{theorem}
For any $p \in (0,1)$ and for $Q \sim Q_2(n,p),$
let the free product factorization of $\pi_1(Q)$ be given by
\[
\pi_1(Q) \cong F * \pi_1(X_1) * \pi_1(X_2) * \cdots * \pi_1(X_\ell),
\]
with $F$ a free group.
With high probability, any $T \in \mathscr{T}_p$ appears as a factor $\pi_1(X_j)$ for some $1 \leq j \leq \ell.$
\end{theorem}
\noindent Our main technical tool will be the following:
\begin{definition}\label{def:hell}
For a cubical sub-complex $T$ of a cubical complex $W\subset Q^n$ denote by $h(T)$ the minimal cubical sub-complex of $W$ so that
\begin{enumerate}
    \item the $1$--skeleton of $h(T)$ is the $1$--skeleton of a $k$--dimensional hypercube
    \item every $2$--face of $W$ that is incident to $T$ is contained in $h(T)$
    \item every $2$--face of $Q^n_2$ with $1$--skeleton in $h(T)$ which is \emph{not} incident to $T$ is in $h(T)$
\end{enumerate}
Also denote by $H(T) \subset W$ as a complete $2$--skeleton of a $k$--dimensional hypercube which is parallel to $h(T)$, so that any 2-face that has an edge in $h(T) \setminus T$ and another edge in $H(T)$ is contained in $W$.
\end{definition}
\noindent We emphasize that $T$ need not be connected in any sense and that $H(T)$ is not unique, but we just need to choose one.
\begin{lemma}\label{lem:hellcube}
For any $p \in (0,1),$ with $Q \sim Q_2(n,p)$ there exists a number $k_p$ so that w.h.p.\ every $((M_p +2)\times k_p)^2$-dimensional cube in $Q^n$ contains fewer than $k_p$ light edges of $Q$.
\end{lemma}
\begin{proof}
We argue by a first moment estimate.  For any $\ell,$ the number of $\ell$--dimensional cubes in $Q^n$ is given by $2^{n-\ell}\binom{n}{\ell}$.  The probability that any such a cube contains $k$ light edges is $O_{k,\ell,p}((1-p)^{nk})$. Hence taking $\ell = ((M_p+2) k)^2,$ if we pick $k$ sufficiently large that $(1-p)^k < \tfrac{1}{2},$ then the expected number of $\ell$--dimensional cubes containing more than $k$ light edges tends to $0$ exponentially in $n.$
\end{proof}

\begin{lemma}\label{lem:heavencube}
Let $p \in (0,1)$ and let $\ell \in \N$ be fixed, then w.h.p.\ for $Q \sim Q_2(n,p),$ every $\ell$-dimensional cube $X$ has a parallel cube $Y$ that has no light edges in $Q$ and for which there are no light edges in $Q$ between $X$ and $Y.$ 
\end{lemma}
\begin{proof}
This is similar to Lemma \ref{lem:hellcube}.  We argue by a first moment estimate.  For any $\ell,$ the number of $\ell$--dimensional cubes in $Q^n$ is given by $2^{n-\ell}\binom{n}{\ell}$.  For a fixed $\ell$--dimensional cube $X \subset Q^n,$ the probability that every parallel $\ell$--dimensional cube $Y$ either
\begin{enumerate}[(i)]
\item contains at least one light edge or
\item contains the endpoint of a light edge between $X$ and $Y$ 
\end{enumerate} 
is at most
\[
((\ell+2)2^{\ell-1})^{n-\ell}
(1-p)^{n(n-\ell)} = o(n^\ell 2^{-n}).
\]
Hence from a first moment estimate, for any fixed $\ell$ and for any $p \in (0,1),$ w.h.p.\ every $\ell$--dimensional cube $X$ has a parallel cube $Y$ that contains no light edges of $Q$ and shares no endpoint of a light edge between $X$ and $Y.$   
\end{proof}

\begin{theorem}\label{thm:heavenandhellexist}
Let $p \in (0,1)$ and $k_p$ as in Lemma \ref{lem:hellcube}. Let $Q \sim Q_2(n,p).$
Let $\overline{Q}$ be $Q$ with all the 2-faces bounded by 4-cycles having no light edges.
With high probability, there are disjoint cubical complexes $\{ \tau_1, \tau_2, \dots, \tau_\ell \}$ in $Q$ so that
\begin{enumerate}[(i)]
    \item the union of $1$-faces over all $\{\tau_j : 1\leq j \leq \ell \}$ is the set of all light $1$-faces,
    \item for each $1 \leq j \leq \ell,$ both $h(\tau_j)$ and $H(\tau_j)$ exist in $\overline{Q},$
    \item for each $1 \leq i \neq j \leq \ell,$ the Hamming distance between the $0$-skeleta of $h(\tau_j)$ and $h(\tau_i)$ is at least 2.
\end{enumerate}
\end{theorem}

\noindent We need the next definition for proving Theorem \ref{thm:heavenandhellexist}. 
\begin{definition}\label{square}
Let $X$ and $Y$ be two subcomplexes of $Q^n$. Define $X \square Y$ to be the face of smallest dimension $Q^m$ such that $Q^m \subset Q^n$ and $X \cup Y \subset Q^m$. Observe that $m\leq n$. 
More in general, let $X_1, X_2, ..., X_l$ be any finite collection of subcomplexes of $Q^n$ and $I=\{1,2, ..., l\}$. We define
\[\square_{i \in I} X_i\]
as the face of smallest dimension $Q^m$ such that $Q^m \subset Q^n$ and such that \[[X_1 \cup X_2 ... \cup X_l] \subset Q^m.\] In this case, $m\leq n$ as well.
\end{definition}

\begin{proof}[Proof of Theorem \ref{thm:heavenandhellexist}]
We first show that for every light $1$-face $e$ there is a cubical complex $\sigma_e$ containing $e$ and having all its $1$-faces light so that $h(\sigma_e)$ exists in $\overline{Q}$. We will then merge these $h(\sigma_e)$ to form the partition claimed to exist in the theorem.

Let $e_1 \coloneqq e$ be any light edge of $Q$.  Let $T_1$ be the cube complex which is the down closure of $e_1.$ Let $X_1$ be the smallest induced complex in $\overline{Q}$  which contains $T_1,$ which contains all $2$-faces of $\overline{Q}$ incident to $e_1$ and whose 1-skeleton is a hypercube.  If $X_1=h(T_1)$, we are finished.  Otherwise, by definition, there must be a $2$-face $f$ of $Q^n_2$ with $1$-skeleton in $X_1$ but which is not itself in $X_1.$  Then, there must be at least one light edge $e_2 \in X_1 \setminus T_1.$

We then define $T_2$ as the induced subcomplex of $Q$ on edges $e_1,e_2$.  Let $X_2$ be the smallest induced complex in $\overline{Q}$ which contains $T_2,$ which contains all $2$-faces of $\overline{Q}$ incident to $T_2$ and whose 1-skeleton is a hypercube.  Once more, if $X_2 = h(T_2),$ we are done.  Otherwise, we proceed inductively by the same argument.

This produces a nested sequence of complexes $\{T_k\}$ each having $k$ edges.  It also produces a sequence of complexes $\{ X_k \}$ such that each $X_k \supset T_k,$ each $X_k$ contains at least $k$ light edges, and such that $X_k$ has the $1$-skeleton of a hypercube of dimension at most $k \times M_p$. By Lemma \ref{lem:hellcube}, with high probability, this sequence must terminate at some $k^* \leq k_p.$  The complex $X_{k^*}=h(T_{k^*})$ by definition, and we define $\sigma_e = T_{k^*}$.

We define a graph $G$ with vertex set given by the collection of $\sigma_e$. Two vertices $\sigma_{e_1}$, $\sigma_{e_2}$ in this graph are connected if the hamming distance between $h(\sigma_{e_1})$ and $h(\sigma_{e_2})$ is less than two. Let $\{ \tau_1, \tau_2, \dots, \tau_\ell \}$ be the unions of the connected components in $G$. Then for each $1 \leq j \leq \ell,$ we construct the hypercube 
\[
\Sigma_{j} = \square_{e \in \tau_j} h(\sigma_{e}). 
\]
It is easy to see that $\Sigma_j=h(\tau_j)$, which implies that $h(\tau_j)$  exists and is exactly $\Sigma_j$. 

The dimension of $\Sigma_j$ is at most 
\begin{equation} 
\sum_{\sigma_e} [\dim(h(\sigma_{e}))+2],
\end{equation}
where the sum is over all $\sigma_{e}$ contained in $\tau_j.$ Therefore by Lemma \ref{lem:hellcube}, each $\tau_j$ has at most $k_p$ edges.
Hence by Lemma \ref{lem:heavencube}, each $H(\tau_j)$ exists as well.
\end{proof}

 \begin{lemma}\label{lem:fromheaven}
 Let $W$ be a subcomplex of $Q_2^n$ and $T$ a subcomplex of $W$.  Suppose that $h(T)$ and $H(T)$ exist in $W$.  
 Let $\hat{W}$ be the complex formed by adding to $W$ the complete $1$-skeleton of $h(T) \square H(T)$ and any $2$-face of $Q_2^n$ with $1$-skeleton in $h(T) \square H(T)$.  Then
 \[
\pi_1(W) 
\cong \pi_1(\hat W) * \pi_1(W \cap (h(T) \square H(T))).
\]
 \end{lemma}
\noindent Recall that for a disconnected cube complex $X,$ we define $\pi(X)$ as the free product of its connected components.
\begin{proof}
Let $\hat{T}$ be all the 1-faces in $T$ and any 2-face of $W$ incident to $T$.
Let $\tilde T$ be the down closure of $\hat T$. Let $S=W \cap (h(T) \square H(T)).$ 

Let $X = (W \setminus \hat{T}) \cap S.$
We claim that $\pi_1(X)\cong 1$. The $1$-faces of $X$ that are in $h(T)$ are not in $T$. Therefore, by Definition \ref{def:hell}, for every edge $e \in X \cap h(T),$ the unique 4-cycle connecting $e$ to $H(T)$ is the boundary of a $2$-face in $X$. Hence, every closed curve in $X$ is homotopic to a curve in $H(T)$.  Since $\pi_1(H(T))\cong 1,$ it follows that $\pi_1(X) \cong 1.$
Therefore, the Siefert-van Kampen theorem states that
\[
\pi_1(W) \cong \pi_1(W \setminus \hat{T}) * \pi_1(S).
\]

We now show that $\pi_1(W \setminus \hat{T}) \cong \pi_1(\hat W)$.  Define a complex $S^*$ as the down closure of all $2$-faces in $Q_2^n$ incident to $T$ with $1$-skeleton in $h(T) \square H(T),$ union with $H(T).$
Any $1$-face $e$ of $S^* \cap (W \setminus \hat T)$ that is in $h(T)$ must be in $h(T) \setminus T$. In particular, there is a $2$-face $f$ containing $e$ which has a $1$-face in $H(T)$.  Hence, any closed curve in $S^* \cap (W \setminus \hat T)$ is homotopic to one in $H(T),$ which is simply connected.  Therefore by the Siefert-van Kampen theorem,
\[
\pi_1(\hat W)
=
\pi_1(S^* \cup (W \setminus \hat{T}))
\cong \pi_1(W \setminus \hat{T}) * \pi_1(S^*).
\]

It remains to evaluate the fundamental group of $S^*$. Any edge in $S^* \cap h(T)$ has a 4-cycle that has an edge in $H(T)$. By construction we know that this 4-cycle has a 2-face added. Therefore any closed curve in $S^*$ is homotopic to a closed curve in $H(T)$. Thus $S^*$ is simply connected because $H(T)$ is by definition.

\end{proof} 

\begin{lemma}\label{lem:heavenhellpi}
 Let $W$ be a subcomplex of $Q_2^n$ and $T$ a subcomplex of $W$.  Suppose that $h(T)$ and $H(T)$ exist in $W$. Let $P_1, ... ,P_m$ be all the pure 2-dimensional strongly connected components completely contained in $T$, such that any 2-face adjacent to the 1-skeleton of any $P_i$ is also contained in $P_i$. 
 Suppose that $T = \cup_{i=1}^m P_i.$
 Then there is a free group $F$ so that 
 \[
\pi_1(W \cap (h(T) \square H(T)))
\cong \pi_1(P_1) * \pi_1(P_2)* \cdots \pi_1(P_m) * F .
\]
 \end{lemma}
 \begin{proof}
 Let $S=W \cap (h(T) \square H(T)).$ Suppose we \textit{fill} $H(T)$ by taking the flag cubical complex of $H(T)$.  The fundamental group of $S$ is unchanged and we can contract $H(T)$ to a point $x$. We denote this complex by $\hat S$. If $e$ is an edge in $T$ then $e$ forms an unfilled triangle with $x$  in $\hat S.$  Let $T_x \subset \hat S$ be the union of $T,$ $x$ and all the edges between $T$ and $x.$ Any edge $f \in h(T)$ which is not contained in $T$ is the base of a filled triangle with $x$ in $\hat S,$ and so any closed curve in $\hat S$ is homotopic to a closed curve in $T_x.$  Hence 
 \[
  \pi_1(S) = \pi_1( \hat S) = \pi_1(T_x)
  = \pi_1(P_1) * \pi_1(P_2)* \cdots \pi_1(P_m) * F
 \]
 where $F$ is a free group.

 \end{proof}

 \begin{theorem}\label{thm:rawstructuretheorem}
 Fix a $p\in(0,1)$. For a $Q \sim Q_2(n,p),$ w.h.p.\ if $\tau_1, \tau_2,...,\tau_\ell$ are as constructed in Theorem \ref{thm:heavenandhellexist}, then with $S_j = C \cap (h(\tau_j) \square H(\tau_j))$ for all $1 \leq j \leq \ell,$
 \[
\pi_1(Q)
\cong  
\pi_1(S_1) * 
\pi_1(S_2) * 
\cdots *
\pi_1(S_\ell).
\]
\end{theorem}
\begin{proof}
Let $\overline{Q}$ be $Q$ with all the 2-faces bounded by 4-cycles having no light edges.  By Lemma \ref{lem:stage2}, all $4$-cycles with no light edges are in $V_3^n$ w.h.p., and so $\pi_1(\overline{Q}) = \pi_1(Q).$
We apply Lemma \ref{lem:fromheaven} inductively to each of the complexes $\tau_j.$ As a result, we have that
\[
\pi_1(\overline{Q})
=
\pi_1(J) *
\pi_1(S_1) * 
\pi_1(S_2) * 
\cdots *
\pi_1(S_\ell),
\]
where $J$ is the complex $\overline{Q}$ together with all $2$-faces in $Q_2^n$ having $1$-skeleton contained in some $h(\tau_j) \square H(\tau_j)$ for some $1 \leq j \leq \ell.$

It just remains to prove that $\pi_1(J) \cong 1.$  The $1$-skeleton of $J$ is $Q^n_1,$ and so it suffices to show that every $4$-cycle in $J$ is contractible.  The only $4$-cycles $x$ in $J$ that do not bound a $2$-face are those that contain a $1$-face $e$ of some $\tau_j$ for $1\leq j \leq \ell$ but which were not contained in $h(\tau_j)\square H(\tau_j).$  However, as $e$ is in $h(\tau_j),$ it has a parallel $1$-face $f$ in $H(\tau_j).$ The unique cube $c = x \square e$ that contains $x$ and $e$ has all 2-faces except for the face bounded by $x$, which implies that $x$ is contractible. 
\end{proof}

\begin{proof}[Proof of Theorem \ref{thm:structure}]
Let $T \in \mathscr{T}_p$ be fixed.  By assumption, there is a $k$--dimensional cube $X$ so that $T$ is a subcomplex of $X$, and we may choose $k$ minimal.  We do not take the full $2$-skeleton for $X,$ but instead we choose exactly those $2$-faces which are either in $T$ or share no edge with $T.$. Note that this makes $X = h(T).$ Let $\phi$ be a cubical embedding of the $2$-skeleton of $X$ into $Q^n_2$.  Define the event $\mathcal{E}_\phi,$ for $Q \sim Q_2(n,p):$
\begin{enumerate}
    \item The $2$-faces of $Q$ that are contained in the $1$-skeleton of $\phi(X)$ are exactly the $2$-faces of $\phi(X).$
    \item No other $2$-face in $Q$ contains a $1$-face of $\phi(T).$
    \item There are no light $1$-faces in $\phi(X\setminus T)$ and no light $1$-faces within Hamming distance $2k_p +2$ of $\phi(T),$ except possibly those in $\phi(T).$ Here, $k_p$ is defined as in Lemma \ref{lem:hellcube}.
\end{enumerate}

We now estimate the probability of $\mathcal{E}_\phi$ under the law of $Q_2(n,p).$ Note that this probability does not depend on $\phi,$ and so these estimates will be uniform in $\phi$.  First, observe that each edge of $\phi(T)$ has degree bounded independently of $n$ on this event, so there are $(e(T)*d)-O(1)$ $2$-faces which must be absent for $\mathcal{E}_\phi$ to hold.  There are $O(1)$ $2$-faces that must be present for $\mathcal{E}_\phi$ to hold, also.  There are also $O(d^{2k_p+2})$ $1$-faces which are contained in the Hamming distance $(2k_p+2)$-neighborhood of $\phi(X)$, which we require to be not light.  As the probability that a $1$-face is light is $O( (1-p)^{n}),$ we conclude that 
\[
\Pr(\mathcal{E}_\phi) = \Theta( (1-p)^{e(T)d} ) = \Omega( 2^{-(1-\epsilon)d}),
\]
for some $\epsilon > 0,$ where the second equality follows from Definition \ref{def:T_p}.
So the expected number of occurrences of $\mathcal{E}_\phi$ goes to infinity exponentially fast as $n\to\infty.$

We can now show that some $\mathcal{E}_\phi$ now occurs with high probability by using a second moment computation (see Corollary 4.3.5 of \cite{AS16}).  Observe that if the Hamming distance of $\phi(X)$ to $\psi(X)$ is greater than $4,$ then the events $\mathcal{E}_{\phi}$ and $\mathcal{E}_{\psi}$ are independent. Let $\psi \sim \phi$ if $\mathcal{E}_{\phi}$ and $\mathcal{E}_{\psi}$ are not independent. Then, 
\[
\Delta_\phi^* = \sum_{\psi \sim \phi} \Pr[ \mathcal{E}_{\psi} | \mathcal{E}_{\phi} ]
\leq \sum_{\psi \sim \phi} 1 = O( d^{O(1)}),
\]
which is much smaller than the expected number of $\mathcal{E}_{\phi}$ that occur (which grows exponentially in $n$).

Hence, with the factorization given by Theorem \ref{thm:rawstructuretheorem},
\[
\pi_1(Q) \cong \pi_1(S_1) * \pi_1(S_2) * \cdots *\pi(S_\ell),
\]
where $S_j = Q \cap (h(\tau_j) \square H(\tau_j))$ and where $\tau_j$ are the complexes from Theorem \ref{thm:heavenandhellexist}.
For any embedding $\phi,$ if $\mathcal{E}_\phi$ occurs, then $\phi(X) \in Q$ is such that $\phi(X)=h(\phi(T))$. Moreover, $\phi(T) = \tau_j$ for some $j$ with $1 \leq j \leq \ell$ as the Hamming distance of $\phi(T)$ to any other light $1$-face is at least $2k_p + 2$. By Lemma \ref{lem:heavenhellpi}, $\pi_1(S_j) \cong F * \pi_1(T)$ for some free group $F.$
\end{proof} 
 
\section{Below the threshold for maximal edges} \label{sec:torsion}

In this section, we discuss in slightly more detail the idea that everything that can happen will happen.

First, we show that in terms of finitely presented groups, everything can happen. The following is well known.

\begin{theorem} \label{thm:simp-to-cube}
Let $S$ be a finite simplicial complex with $k$ vertices. Then $S$ is homeomorphic to a cube complex. Indeed, $S$ is homeomorphic to a  subcomplex of the $k$-dimensional cube. 
\end{theorem}

This is well known. See, for example, Appendix A of Davis's book \cite{Davis08}, or the discussion in Section 2 of Babson, Billera, and Chen's paper \cite{BBC97}.

Theorem \ref{thm:simp-to-cube} has a the following immediate corollary.

\begin{corollary}
Let $G$ be any finitely-presented group. Then there exists a number $P_G >0$ such that whenever $0 < p < P_G$ and $Q \sim Q_2(n,p)$, we have that $G$ exists as a free factor in $\pi_1( Q)$ with high probability.
\end{corollary}

\begin{figure}
    \centering
    \includegraphics[scale = .4]{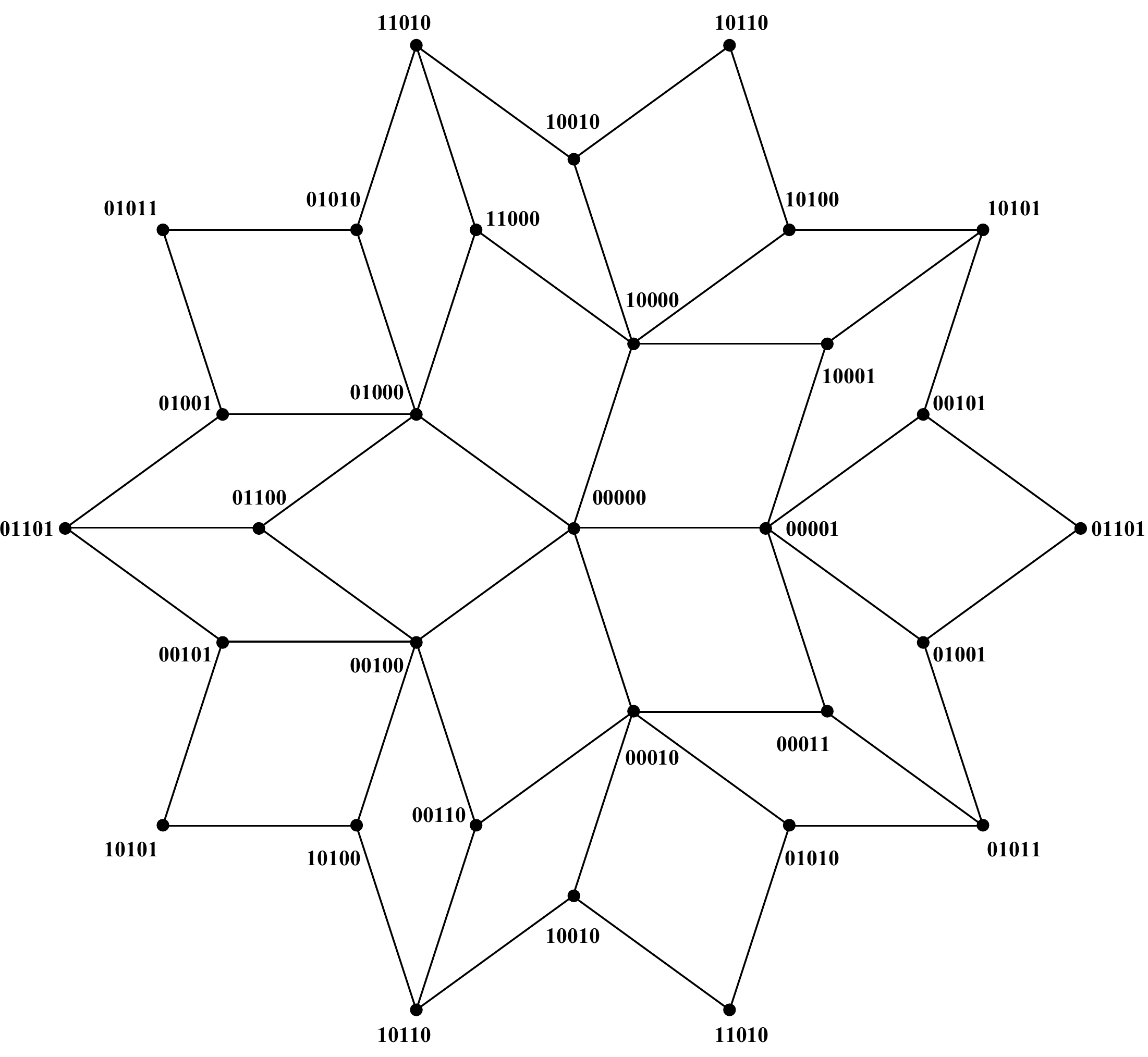}
    \caption{Diagram of $T_2^\square$ with the vertices labeled.}
    \label{F:isjustbeutiful}
\end{figure}

We give concrete bounds on $P_G$ for a few groups $G$ in the following. The constructions used in Theorem \ref{thm:12} and Theorem \ref{thm:20} and Figure \ref{F:isjustbeutiful} are from joint work of Dejan Govc and the third author of this paper. 

\begin{theorem}\label{thm:13}
Let $T^2$ be the $2$-dimensional torus. For $0 < p < (1-(1/2)^{1/2^5})\approx 0.021428,$ $\pi_1(T_2) \cong \Z \times \Z$ is a free factor of $\pi_1(Q)$ for $Q \sim Q_2(n,p)$ with high probability. 
\end{theorem} 

\begin{proof}

We will create a cubical complex which is a subcomplex of $Q_2^4$ and which is homeomorphic to $T^2.$  
Let $T_1^\square$ be the minimal cubical subcomplex of $Q_2^4$ with 2-faces given by
\[
 \begin{aligned}
\{&(*,0,0,*), (0, *, 0, *), (*, 1, 0, *), (1, *, 0, *),\\&(*, 0, *, 1), (0, *, *, 1), (*, 1, *, 1), (1, *, *, 1),\\& (*, 0, 1, *), (0, *, 1, *), (*, 1, 1, *), (1, *, 1, *),\\& (*, 0, *, 0), (0, *, *, 0), (*, 1, *, 0), (1, *, *, 0)\}.
\end{aligned}
\]
Observe that $T_1^\square$ has $e(T_1^\square) =32$. 
Hence, from Theorem \ref{thm:structure}, for $0<p \neq 0$, if $ p < (1-(1/2)^{1/32}),$ the fundamental group of $Q \sim Q_2(n,p)$ has a copy of $\Z \times \Z$ in its free product factorization with high probability.

\end{proof} 

\begin{theorem}\label{thm:12}
Let $T_2$ be the projective plane. For $p \neq 0$,  if $ p < (1-(1/2)^{1/40})\approx 0.017179,$ $\pi_1(T_2) \cong \Z / 2\Z$ is a free factor of $\pi_1(Q)$ for $Q \sim Q_2(n,p)$ with high probability.
\end{theorem}

\begin{proof}
Let $T_2^\square$ be the minimal cubical subcomplex of $Q_2^5$ with 2-faces given by
\[
 \begin{aligned}
\{&(0,0, 0, *, *), (0, 0, 1, *, *), (0, 0, *, 1, *), (0, 0, *, *, 1),\\ &(0, 1, *, *, 0), (0, *, 0, 0, *), (0, *, 1, *, 0), (0, *, *, 0, 0),\\ &(0, *, *, 1, 0), (1, 0, *, 0, *), (1, *, 0, 0, *), (1, *, 0, *, 0),\\ & (*, 0, 0, *, 0), (*, 0, 1, 0, *), (*, 0, *, 0, 0), (*, 0, *, 0, 1),\\ & (*, 1, 0, 0, *), (*, *, 0, 0, 1), (*, 1, 0, *, 0), (*, *, 0, 1, 0)\},
\end{aligned}
\]
so that $e(T_2^\square) =40$.  
Hence, for $p \neq 0$, if $ p < (1-(1/2)^{1/40})\approx 0.017179,$ the fundamental group of $Q \sim Q_2(n,p)$ has a torsion group in its free product factorization with high probability.
\end{proof}

\begin{theorem}\label{thm:20}
Let $G = < x, y \mid x y x^{-1} y>$ be the fundamental group of the Klein bottle. If $0 < p < (1-(1/2)^{1/56})\approx 0.01230134$ and $Q \sim Q_2(n,p)$, then with high probability \ $\pi_1(Q)$ has $G$ as a free factor.
\end{theorem}

\begin{proof}
Let $T_3^\square$ be the minimal cubical subcomplex of $Q_2^5$ with 2-faces given by
\[
 \begin{aligned}
\{&(*,*,1,0,0), (*,*,0,0,0), (0,*,*,0,0), (1,*,*,0,0),\\&(0,1,*,*,0),(1,1,*,*,0),(*,1,0,*,0),(*,1,1,*,0),\\ &(0,*,*,1,0),(1,*,*,1,0),(*,*,1,1,0),(*,0,*,1,0),\\ &(0,*,0,1,*),(1,*,0,1,*),(*,0,0,1,*),(*,1,0,1,*), \\ &(*,1,*,1,1),(0,*,*,1,1),(1,*,*,1,1),(*,*,1,1,1),\\ &(0,0,*,*,1),(1,0,*,*,1),(*,0,0,*,1),(*,0,1,*,1),\\ &(0,0,*,0,*),(1,0,*,0,*),(*,0,0,0,*),(*,0,1,0,*)\}.
\end{aligned}
\]
Then $e(T_3^\square) = 56$. 
\end{proof}

\subsection*{Acknowledgements}
We would like to thank Yuval Peled for interesting conversations that helped launch this work. We thank Dejan Govc for the constructions used in Theorem \ref{thm:12} and Theorem \ref{thm:20} in Section 5, and for Figure \ref{F:isjustbeutiful}.

\bibliographystyle{abbrvnat}
\bibliography{bib}

\end{document}